\documentclass[12pt]{article}
\usepackage{amsmath,amsthm,amssymb,amsfonts}
\usepackage[numbers,sort&compress]{natbib}
\usepackage{color,colordvi}
\usepackage{soul}
\usepackage{fullpage}

\usepackage{amsmath,amsfonts,amsthm,mathtools,bm}
\usepackage{graphicx}
\usepackage[colorlinks=true, allcolors=blue]{hyperref}
\usepackage[capitalise]{cleveref}
\usepackage{standalone}
\usepackage{indentfirst}
\usepackage[affil-it]{authblk}
\usepackage{diagbox}
\usepackage{changes}
\usepackage{enumitem}
\usepackage{bbm}
\usepackage{orcidlink}

\newtheorem{theorem}{Theorem}[section]
\newtheorem{lemma}[theorem]{Lemma}
\newtheorem{corollary}[theorem]{Corollary}
\newtheorem{proposition}[theorem]{Proposition}
\newtheorem{claim}[theorem]{Claim}

\newtheorem{conjecture}[theorem]{Conjecture}
\newtheorem{observation}[theorem]{Observation}

\newtheorem*{claim*}{Claim}

\newcommand{\textotherwise}{\text{otherwise}}

\newcommand{\QQ}{\mathbb{Q}}
\newcommand{\ZZ}{\mathbb{Z}}

\newcommand{\cG}{\mathcal{G}}

\newcommand{\cQ}{\mathcal{Q}}

\DeclareMathOperator{\diag}{diag}

\DeclareMathOperator{\CAN}{CAN}
\DeclareMathOperator{\sCAN}{semi-CAN}
\DeclareMathOperator{\FRI}{FRI}
\DeclareMathOperator{\FCI}{FCI}
\DeclareMathOperator{\IRI}{IRI}
\DeclareMathOperator{\ICI}{ICI}

\DeclarePairedDelimiter{\card}{\lvert}{\rvert}
\DeclarePairedDelimiter{\set}{\lbrace}{\rbrace}

\DeclarePairedDelimiter{\paren}{\lparen}{\rparen}

\DeclarePairedDelimiter{\floor}{\lfloor}{\rfloor}
\DeclarePairedDelimiter{\ceil}{\lceil}{\rceil}

\title{Almost all graphs have no cospectral mate with fixed level}

\author[a]{Wei Wang}
\affil[a]{School of Mathematics and Statistics, Xi'an Jiaotong University, Xi'an, China.
\href{mailto:wang\_weiw@xjtu.edu.cn}{wang\_weiw@xjtu.edu.cn}}
\author[b]{Da Zhao~\orcidlink{0000-0002-9582-0778}}
\affil[b]{School of Mathematics, East China University of Science and Technology, Shanghai, China. \href{mailto:zhaoda@ecust.edu.cn}{zhaoda@ecust.edu.cn}}

\date{}

\begin{document}
\maketitle

\begin{abstract} 
    Haemers conjectures that almost all graphs are determined by their spectra. 
    Suppose $G \sim \mathcal{G}(n, p)$ is a random graph with each edge chosen independently with probability $p$ with $0 < p < 1$.
    Then
    $$\Pr(G \text{ is not controllable}) + \sum_{\ell = 2}^{n^{n^2}} \Pr(G \text{ has a generalized copsectral mate with level } \ell) \to 0$$ as $n \to \infty$ implies that almost all graphs are determined by their generalized spectra. 
    It is known that almost all graphs are controllable. 
    We show that almost all graphs have no cospectral mate with fixed level $\ell$, namely $$\Pr(G \text{ has a copsectral mate with level } \ell) \to 0$$ as $n \to \infty$ for every $\ell \geq 2$. 
    Consequently, 
    $$\Pr(G \text{ has a generalized copsectral mate with level } \ell) \to 0$$ as $n \to \infty$ for every $\ell \geq 2$. 
    The result can also be interpreted in the framework of random integral matrices. 
\end{abstract}



\section{Introduction}

The graph isomorphism problem asks whether two given graphs are isomorphic or not. 
The complexity-theoretic results suggest that the isomorphism problem is not NP-complete~\cite{boppana1987DoesCoNPHave, schoning1988GraphIsomorphismLow}. 
On the other hand, no polynomial-time algorithm is known so far. 
However, in practice, the graph isomorphism problem can be solved efficiently for almost all graphs. 
Babai--Ku\v{c}era~\cite{babai1979CanonicalLabellingGraphs} and Babai--Erd\H{o}s--Seklow~\cite{babai1980RandomGraphIsomorphism} show that almost all graphs are determined by their color refinements. 
An improvement of the color refinement algorithm is the $k$-dimensional Weisfeiler--Lehman ($k$-WL) algorithm whose time complexity is $\mathcal{O}(n^{k+1} \log n)$. 
This provides combinatorial reasons for why the graph isomorphism problem is solved efficiently in practice. 
The matrices of a graph are natural algebraic tools to study a graph. 
The distribution of eigenvalues of a random matrix is well-known, say Erd\H{o}s--Yau~\cite{erdos2017dynamical} and Tao~\cite{tao2023topics}. 
Haemers conjectures that almost all graphs are determined by their eigenvalues. 
In other words, almost all graphs have no cospectral mate. 
This may serve as algebraic reasons for why the graph ismorphism problem is solved efficiently in practice. 
In this paper, we prove that almost all graphs have no cospectral mate with a fixed level.

A \emph{graph} is a tuple $G = (V, E)$, where $V$ is the \emph{vertex} set and $E \subseteq \binom{V}{2}$ is the \emph{edge} set.
The \emph{adjacency matrix} of $G$ is defined by $A = A_G = (a_{i, j})_{i, j \in V}$, where
\begin{align}
    a_{i,j} = 
    \begin{cases}
        1, & i \sim j \text{ is an edge}, \\
        0, & \textotherwise.
    \end{cases}
\end{align}
The \emph{spectrum} of a graph $G$ is the multiset of all eigenvalues of $A_G$. 

Let $S_n(\ZZ)$ be the set of integral symmetric matrices of order $n$;
let $M_n(\ZZ)$ be the set of integral matrices of order $n$; 
let $O_n(\QQ)$ be the set of rational orthogonal matrix of order $n$. Given a rational matrix $M$, the \emph{level} $\ell(M)$ of $M$ is the smallest positive integer $\ell$ such that $\ell M$ is an integer matrix. 
For every positive integer $\ell$, let $O_{n, \ell}(\QQ)$ be the set of rational orthogonal matrices with level $\ell$. 
Note that the rational orthogonal matrices with level $1$ are exactly the signed permutation matrices. 

Two graphs $G_1 = (V_1, E_1)$ and $G_2 = (V_2, E_2)$ are called \emph{isomorphic} if there exists a bijection $\pi : V_1 \to V_2$ such that $uv \in E_1$ if and only if $\pi(u)\pi(v) \in E_2$. 
Two non-isomorphic graphs are called \emph{cospectral} if their spectra are identical. 
And we say that one graph is a \emph{cospectral mate} of the other. 
In particular, this implies that there exists an orthogonal matrix $Q$ such that $Q^\top A Q = B$, where $A$ and $B$ are the adjacency matrices of the two graphs. 
In general the orthogonal matrix $Q$ may not be rational. 
In the case that $Q$ is a rational orthogonal matrix with level $\ell$, we say that one graph is a cospcetral mate of the other with level $\ell$. 

We consider the random graph model $\cG(n, p)$. 
Namely, $G \sim \cG(n, p)$ is a graph on $n$ vertices with each edge chosen independently with probability $p$. 

\begin{conjecture}[{Haemers conjecture~\cite{vandam2003WhichGraphsAre,haemers2016almost}}]
    Almost all graphs are determined by their spectra. 
    In other words, almost all graphs have no cospectral mate. 
\end{conjecture}

The \emph{complement} of a graph $G = (V, E)$ is the graph $\overline{G} = (V, \overline{E})$, where $\overline{E} = \binom{V}{2} \setminus E$. 
The \emph{generalized spectrum} of a graph $G$ is the multiset of all eigenvalues of $G$ combined with the multiset of all eigenvalues of $\overline{G}$. 
Two non-isomorphic graphs are called generalized cospectral if their generalized spectra are identical. 
An orthogonal matrix is called regular if the sum of the entries in each row is $1$.
The \emph{walk matrix} of a graph $G = (V, E)$ is defined by
\begin{align}
    W(G) = \begin{bmatrix}
        e & Ae & \cdots & A^{n-1}e
    \end{bmatrix},
\end{align}
where $e$ is the all-one column vector and $n = \card{V}$.
A graph $G$ is called \emph{controllable} if $W(G)$ is of full rank. 
The following theorem characterizes generalized cospectral graphs. 

\begin{theorem}[{\cite{johnson1980NoteCospectralGraphs}~\cite{wang2006SufficientConditionFamily}}]
    Let $G$ and $H$ be two non-isomorphic graphs. 
    Let $A$ and $B$ be their adjacency matrices respectively. 
    If $G$ and $H$ are generalized cospectral, then there exists a regular orthogonal matrix $Q$ such that $Q^\top A Q = B$. 
    In particular, if $G$ is controllable, then $Q$ is unique and rational. 
    Moreover, $\ell(Q) \mid d_n(W(G))$, where $d_n(W(G))$ is the last invariant factor in the Smith normal form of $W(G)$.
\end{theorem}

\begin{corollary}
    If
    \begin{align}\label{eq:across}
        \Pr(G \text{ is not controllable}) + \sum_{\ell = 2}^{n^{n^2}} \Pr(G \text{ has a generalized copsectral mate with level } \ell) \to 0
    \end{align}
    as $n \to \infty$, then almost all graphs are determined by their generalized spectra. 
\end{corollary}

The first term is in~\cref{eq:across} is shown approaching $0$ as $n \to \infty$ by O'Rourke--Touri. 

\begin{theorem}[{\cite{orourke2016ConjectureGodsilConcerning}}]
    Almost all graphs are controllable. 
\end{theorem}

Based on the classification of regular rational orthogonal matrix of level $2$, Wang--Xu proved the following theorem. 

\begin{theorem}[{\cite{wang2010AsymptoticBehaviorGraphs}}]
    Let $G \sim \cG(n, 1/2)$. 
    Then
    \begin{align}
        \Pr(G \text{ has a generalized cospectral mate with level } 2) \to 0
    \end{align}
    as $n \to \infty$. 
\end{theorem}

In this paper, we generalize the above theorem in three directions simultaneously. 
The graph model is extended from $\cG(n, 1/2)$ to $\cG(n, p)$ with $0 < p < 1$. 
The level can be an arbitrary positive integer $\ell > 1$. 
And the generalized cospectral condition is relaxed to cospectral condition. 
In particular, we show that every summand in~\cref{eq:across} tends to $0$ as $n \to \infty$. 

\begin{theorem}\label{thm:guide}
    Let $0 < p < 1$ and $G \sim \cG(n, p)$. 
    Let $\ell \geq 2$ be an integer. 
    Then almost all graphs have no cospectral mate with level $\ell'$ with $\ell' \mid \ell$. 
    Namely,
    \begin{align}
        \Pr(G \text{ has a cospectral mate with level } \ell' \text{ such that } \ell' \mid \ell) \to 0
    \end{align}
    as $n \to \infty$. 
    In particular, 
    \begin{align}
        \Pr(G \text{ has a cospectral mate with level } \ell) \to 0
    \end{align}
    as $n \to \infty$. 
    Consequently, almost all graphs have no generalized cospectral mate with level $\ell'$ with $\ell' \mid \ell$. 
    Namely,
    \begin{align}
        \Pr(G \text{ has a generalized cospectral mate with level } \ell) \to 0
    \end{align}
    as $n \to \infty$. 
\end{theorem}

\section{Proof}

The proof is divided into four parts. 
Firstly, we reduce the rational orthogonal matrix to canonical form. 
Secondly, for a fixed $Q$ of canonical form, we estimate the probability that $\Pr(Q^\top A Q \in M_n(\ZZ))$. 
Thirdly, we estimate the number of $Q$ in canonical form. 
Lastly, we combine the estimates to prove the main theorem. 

For $Q \in O_{n, \ell}(\QQ)$, there exists a pair of signed permutations matrices $(P_R, P_C)$ such that $P_R^\top Q P_C$ is of the block diagonal form
\begin{align}
    \begin{bmatrix}
        Q_s & \\
         & I_{n-s}
    \end{bmatrix},
\end{align}
where $Q_s$ is a rational orthogonal matrix of order $s$ and level $\ell$ with entries being zero or fractions.  
We call it a \emph{canonical form} of $Q$. 
Note that the canonical form a rational orthogonal matrix is not unique. 
We denote by $\CAN(n, \ell; s)$ the subset of $O_{n, \ell}(\QQ)$ of block diagonal form $\diag(Q_s, I_{n-s})$, and by $\sCAN(n, \ell; s)$ the subset of $O_{n, \ell}(\QQ)$ of block diagonal form $\diag(Q_s, P_{n-s})$, where $Q_s$ is a rational orthogonal matrix of order $s$ with entries being zero or fractions, and $P_{n-s}$ is a signed permutation matrix of order $n-s$. 

\begin{lemma}\label{lem:bottle}
    Fix $Q \in \CAN(n, \ell; s)$. 
    Let $A$ be the adjacency matrix of a random graph $G \sim \cG(n, p)$. 
    Then
    \begin{align}
        \Pr(Q^\top A Q \in M_n(\ZZ)) \leq \hat{p}^{\frac{s}{2 \ell^4} \cdot (\frac{s}{2\ell^4} + n - s - 1)},
    \end{align}
    where $\hat{p} = \max\set{p, 1-p}$. 
\end{lemma}

Fix $Q \in \CAN(n, \ell; s)$. 
For a column index $j$ of $Q$, we define $K(j)$ to be the set of row indices of non-zero entries of the $j$-th column of $Q$. 
Next we shall choose two subsets of indices $I \subseteq \set{1,2, \ldots, s}$ and $J \subseteq \set{1,2, \ldots, n}$ such that the followings hold.
\begin{enumerate}
    \item It holds $K(i) \cap K(j) = \emptyset$ for $i, j \in I \cup J, i \neq j$. 
    \item The Cartesian products $K(i) \times K(j), i \in I, j \in J$ are disjoint subsets of $\set{1,2, \ldots, s} \times \set{1,2, \ldots, n}$;
    \item The size of $I$ and $J$ are large enough in the sense that $\card{I} = \Omega(s)$ and $\card{J} = \Omega(n)$.
\end{enumerate}

For every index $i \in \set{1,2, \ldots, s}$, we define the set $N(i) = \set{j \in \set{1,2, \ldots, s} \mid K(j) \cap K(i) \neq \emptyset}$. 

\begin{claim}
    For every $i \in \set{1,2, \ldots, s}$, we have $\card{N(i)} \leq \ell^4$. 
\end{claim}

\begin{proof}
    Note that there are at most $\ell^2$ non-zero entries in each row or column of $Q$. 
    Therefore, $\card{K(i)} \leq \ell^2$. 
    For every $k \in K(i)$, there exist at most $\ell^2$ non-zero entries in the $k$-th row of $Q$. 
    Hence, $L(k) \coloneqq \card{\set{j \in \set{1,2, \ldots, s} \mid k \in K(j)}} \leq \ell^2$.
    Since $N(i) = \cup_{k \in K(i)} L(k)$, we get $\card{N(i)} \leq \ell^2 \times \ell^2 = \ell^4$.
\end{proof}

Define $S_1 = \set{1,2, \ldots, s}$. 
We take an arbitrary element $i_1 \in S_1$, and put the index $i_1$ into $I$. 
We define $S_2 = S_1 \setminus N(i_1)$. 
If $S_2 \neq \emptyset$, then we take an arbitrary element $i_2 \in S_2$, and put the index $i_2$ into $J$. 
We define $S_3 = S_2 \setminus N(i_2)$. 
If $S_3 \neq \emptyset$, then we take an arbitrary element $i_3 \in S_2$, and put the index $i_3$ into $I$. 
We repeat the above process until $S_{t+1} = \emptyset$ for some positive integer $t$.
Finally, we put $\set{s+1, s+2, \ldots, n}$ to $J$. 

\begin{claim}
    We have $\card{I} \geq \ceil{\frac{1}{2}\ceil{\frac{s}{\ell^4}}} \geq \frac{s}{2\ell^4}$, and $\card{J} \geq \floor{\frac{1}{2}\ceil{\frac{s}{\ell^4}}} + n-s \geq \frac{s}{2\ell^4} + n - s - 1$.
\end{claim}

\begin{proof}
    Since $\card{N(i)} \leq \ell^4$ for every $i \in \set{1,2, \ldots, s}$, we have $S_t \neq \emptyset$ for $t = \ceil{\frac{s}{\ell^4}}$. 
    Hence, $\card{I} \geq \ceil{\frac{t}{2}}$, and $\card{J} \geq \floor{\frac{t}{2}} + n-s$. 
\end{proof}

Let $q_j, j \in \set{1,2, \ldots, n}$ be the $j$-th column of $Q$. 
Let $q_{i,j}$ be the $i$-th row of $q_i$. 
Let $b_{i,j} = q_i^\top A q_j$ be the $(i,j)$-entry of $Q^\top A Q$. 

\begin{proposition}
    It holds that
    \begin{align}
        \Pr(b_{i,j} \in \ZZ) \leq \hat{p}
    \end{align}
    for $i \in I, j \in J$. 
    And $b_{i,j}$ are independent random variables for $i \in I, j \in J$. 
\end{proposition}

\begin{proof}
    Note that
    \begin{align}
        b_{i, j} &= (Q^\top A Q)_{i, j} \\
            &= \sum_{1 \leq u, v \leq n} q_{u,i} a_{u, v} q_{v, j} \\
            &= \sum_{u \in K(i), v \in K(j)} a_{u, v} q_{u,i} q_{v, j}. 
    \end{align}
    So $b_{i, j}$ depends only on the variables $a_{u, v}$ for $(u, v) \in K(i) \times K(j)$. 
    Since $K(i) \times K(j), i \in I, j \in J$ are disjoint subsets of $\set{1,2, \ldots, s} \times \set{1,2, \ldots, n}$ and $K(i) \cap K(j) = \emptyset$ for $i, j \in I \cup J$ and $i \neq j$, we have that $b_{i,j}$ are independent for $i \in I, j \in J$. 

    Let $X = (x_{u, v})_{u \in K(i), v \in K(j)}$ be a $(0, 1)$-array with indices $K(i) \times K(j)$. 
    We define an involution on all such arrays.  
    Fix $u_0 \in K(i), v_0 \in K(j)$. 
    Define $X' = (x_{u, v}')_{u \in K(i), v \in K(j)}$, where
    \begin{align}
        x_{u, v}' = 
        \begin{cases}
            x_{u, v}, & (u, v) \neq (u_0, v_0), \\
            1-x_{u, v}, & (u, v) = (u_0, v_0).
        \end{cases}
    \end{align}
    Note that $0 < q_{u_0, i} < 1$ and $0 < q_{v_0, j} \leq 1$. 
    So if $b_{i, j} \in \ZZ$ for $A|_{K(i) \times K(j)} = X$, then $b_{i, j} \notin \ZZ$ for $A|_{K(i) \times K(j)} = X'$. 
    Namely,
    \begin{align}
        \Pr(b_{i,j} \in \ZZ \mid a_{u, v} = x_{u, v}, (u, v) \neq (u_0, v_0)) \leq \max(p, 1-p) = \hat{p}. 
    \end{align}
    By total probability theorem we get
    \begin{align}
        \Pr(b_{i,j} \in \ZZ) \leq \hat{p}. 
    \end{align}
\end{proof}

\begin{proof}[{Proof of~\cref{lem:bottle}}]
    For large $n$, we have
    \begin{align}
        \Pr(Q^\top A Q \in M_n(\ZZ))
        &\leq \prod_{i \in I, j \in J} \Pr(b_{i, j} \in \ZZ) \\
        &\leq \hat{p}^{\card{I} \cdot \card{J}} \\
        &\leq \hat{p}^{\frac{s}{2 \ell^4} \cdot (\frac{s}{2\ell^4} + n - s - 1)}.
    \end{align}
\end{proof}

Next we estimate the number of $Q$ in canonical form. 

\begin{lemma}
    Let $\ell$ be a positive integer. 
    Then
    \begin{align}
        \sum_{\ell' \mid \ell} \card{O_{n, \ell'}(\QQ)} \leq (2n)^{\ell^2 n }.
    \end{align}
\end{lemma}

\begin{proof}
    Let $Q \in O_{n, \ell'}(\QQ)$ with $\ell' \mid \ell$. 
    Note that $\ell Q$ is an integral matrix. 
    Moreover, the sum of squares of entries in each column of $\ell Q$ is $\ell^2$. 
    The number of feasible column vector is at most $n^{\ell^2} 2^{\ell^2}$ with signs of entries taking into account. 
    Therefore, $\sum_{\ell' \mid \ell} \card{O_{n, \ell'}(\QQ)} \leq \paren*{(2n)^{\ell^2 }}^n = (2n)^{\ell^2 n}$.
\end{proof}

\begin{corollary}
    \begin{align}
        \sum_{\ell' \mid \ell}\card{\CAN(n, \ell; s)} \leq (2s)^{\ell^2 s}.
    \end{align}
\end{corollary}

For $A \in S_n(\ZZ)$, define 
\begin{align}
    \cQ(A) \coloneqq \set{Q \in O_n(\QQ) \mid Q^\top A Q \in M_n(\ZZ)}.
\end{align}

\begin{observation}
    Let $A \in S_n(\ZZ)$.
    Suppose $Q \in \cQ(A)$. 
    Then $QP \in \cQ(A)$ for every signed permutation matrix $P$ of order $n$. 
\end{observation}

For $Q \in O_n(\QQ)$, we denote by $\FRI(Q)$ ($\FCI(Q)$) the index set of rows (columns) with fractional numbers, and by $\IRI(Q)$ ($\ICI(Q)$) the index set of rows (columns) with only integral numbers (namely $0$ and $\pm 1$).

\begin{observation}
    Let $A \in S_n(\ZZ)$ and $Q \in \cQ(A)$. 
    Suppose $P_R$ and $P_C$ are permutation matrices mapping $\FRI(Q)$ and $\FCI(Q)$ to $\set{1,2, \ldots, s}$ respectively, and $P_R^\top Q P_C \in \sCAN(n, \ell; s)$. 
    Then there exists $Q' \in \cQ(A)$ such that $P_R^\top Q' P_C \in \CAN(n, \ell; s)$. 
\end{observation}

Note that the number of choices for the pair of permutation matrices $P_R$ and $P_C$ is at most $(n (n-1) \cdots (n-s+1))^2 \leq n^{2s}$. 

We are prepared to prove the main theorem.

\begin{proof}[{Proof of~\cref{thm:guide}}]
    For large $n$, we have
    \begin{align}
        &\Pr(G \text{ has a cospectral mate with level } \ell' \text{ such that } \ell' \mid \ell) \\
        \leq&\Pr(\exists Q \in O_{n, \ell'}(\QQ) : 2 \leq \ell' \mid \ell, Q^\top A Q \in M_n(\ZZ)) \\
        \leq& \sum_{s = 2}^n \sum_{\ell' \mid \ell} \sum_{Q \in \CAN(n, \ell'; s)} n^{2s} \Pr(Q^\top A Q \in M_n(\ZZ)) \\
        \leq& \sum_{s = 2}^n n^{2s} (2s)^{\ell^2 s} \hat{p}^{\frac{s}{2 \ell^4} \cdot (\frac{s}{2\ell^4} + n - s - 1)} \\
        \leq& \sum_{s = 2}^n n^{2s} (2n)^{\ell^2 s} \hat{p}^{\frac{s}{2 \ell^4} \cdot \frac{n-1}{2\ell^4}} \\
        =& \sum_{s = 2}^n \epsilon_n^s \leq \frac{\epsilon_n^2}{1 - \epsilon_n} \to 0, \text{ as } n \to \infty,
    \end{align}
    where $\epsilon_n = n^{2} (2n)^{\ell^2} \hat{p}^{\frac{1}{2 \ell^4 \cdot 2 \ell^4} (n-1)}$.
\end{proof}

\section*{Acknowledgements}

Wei WANG is supported in part by the National Key Research and Development Program of China 2023YFA1010203, and the National Natural Science Foundation of China (No. 12371357). 
Da ZHAO is supported in part by the National Natural Science Foundation of China (No. 12471324, No. 12501459, No. 12571353), and the Natural Science Foundation of Shanghai, Shanghai Sailing Program (No. 24YF2709000). 

\bibliographystyle{alpha}
\bibliography{ref}

\end{document}